# From Integer Sequences to Block Designs via Counting Walks in Graphs

Ernesto Estrada[1,2] and José A. de la Peña[2]

[1]Department of Mathematics and Statistics, University of Strathclyde, Glasgow G1 1XH, U.K., [2]Centro de Investigación en Matemáticas (CIMAT), A. C., Guanajuato 36240, México.

## ABSTRACT

We define numbers of the type $O_j(N) = N^0 - N^1 + N^2 - \cdots + N^{2j}$ and $E_j(N) = -N^0 + N^1 - N^2 + \cdots + N^{2j+1}$ ($j = 0, 1, 2, \cdots$) and the corresponding integer sequences. We prove that these integer sequences, e.g., $S_O(N) = O_0(N), O_1(N), \cdots, O_r(N), \cdots$ and $S_E(N) = E_0(N), E_1(N), \cdots, E_r(N), \cdots$ correspond to the number of odd and even walks in complete graphs $K_N$. We then prove that there is a unique family of graphs which have exactly the same sequence of odd walks between connected nodes and of even walks between pairs of nodes at distance two, respectively. These graphs are obtained as the Kronecker product $G_{2n} = K_2 \otimes K_n$. We show that they are the incidence graphs of block designs, are distance-regular and Ramanujan graphs.

## 1. Introduction

Integer sequences arise from a few different sources, such as enumeration problems, number theory, game theory, physics, and so forth [1]. Enumeration problems in graph theory are an immense source of integer sequences as can be seen from an inspection of the *Online Encyclopedia of Integer Sequences* [2, 3]. In that particular case an investigation in the context of graph theory or combinatorics gives rise to a sequence of integers, which is then properly investigated. However, here we pose a different problem. For instance, let $N \in \mathbb{N}$ be natural number ($N > 1$) and let us define the following numbers ($j = 0, 1, 2, \cdots$):

$$O_j(N) = N^0 - N^1 + N^2 - \cdots + N^{2j}, \tag{1}$$

$$E_j(N) = -N^0 + N^1 - N^2 + \cdots + N^{2j+1}. \tag{2}$$

Using $O_j(N)$ and $E_j(N)$ let us now introduce the following sequences:

$$S_O(N) = O_0(N), O_1(N), \cdots, O_r(N), \cdots \tag{3}$$

$$S_E(N) = E_0(N), E_1(N), \cdots, E_r(N), \cdots. \tag{4}$$

There are two $S_O(N)$ and three $S_E(N)$ sequences reported in the *Online Encyclopedia of Integer Sequences* (OEIS) [3]. They are:

| Sequence | OEIS code | Current work |
|---|---|---|
| 1, 3, 11, 43, 171,... | A007583 | $S_O(2)$ |
| 1, 7, 61, 547, 4921, ... | A066443 | $S_O(3)$ |
| 0, 1, 5, 21, 85, 341... | A002450 | $0, S_E(2)$ [a] |
| 0, 2, 20, 182, 1640,... | A125857 | $0, S_E(3)$ [a] |
| 0, 3, 51, 819, 13107,... | A182512 | $0, S_E(4)$ [a] |

[a]The even sequences here do not include zero as in the OIES ones.

They appear, however, related to different mathematical objects. For instance, A007583, A002450 and A182512 appear associated to the formulae $(2^{2n+1}+1)/3$, $(4^n-1)/3$ and $(16^n-1)/5$, respectively. However, A066443 appears as the number of distinct walks of length $2n+1$ along edges of a unit cube between two fixed adjacent vertices and A125857 represents the numbers whose base 9 representation is 22222222.......2. For $N>3$ no sequence of the type $S_o(N)$ is reported in the OEIS and the same is true for $S_E(N)$ when $N>4$. So, the questions are: Is there a general framework in which all of these sequences can be grouped together? Are these sequences related to any property of a certain kind of graph?

Of course, these questions are not always possible to be answered. However, in the particular case that we can associate such sequences with a property of a type of graph we can yet ask another fundamental question. Given an integer sequence that represents a property of a given kind of graph, can we construct another family of graphs having the same sequence for this property? These questions are investigated here for the particular sequences (3) and (4). We find here that they correspond to the number of odd and even-length walks between pairs of nodes in complete graphs. We then prove that for each complete graph $K_n$ there is a unique graph, not isomorphic to $K_n$, which has exactly the same sequence of walks. Further investigations of these graphs allow us to prove that they are the incidence graphs of certain symmetric block designs [4-6]. In addition, the graphs constructed in this way are distance-regular [7] and Ramanujan graphs [8], showing potentialities in many different areas of applications. In short, we start here with some integer sequences and end up with some graphs that are related to symmetric block designs with several interesting properties. The connection between these areas is provided by the sequences of walks in complete graphs. In this work all the definitions are given the first time that the concept is found in the text.

## 2. Integer Sequences and Walks in Complete Graphs

We start by finding the general expressions for the numbers $O_j(N)$ and $E_j(N)$.

**Lemma 1:** The numbers $O_j(N)$ and $E_j(N)$ are positive integers given by the following formulae:

$$O_j(N) = \frac{N^{2j+1} + 1}{N+1}, \tag{5}$$

$$E_j(N) = \frac{N^{2j+2} - 1}{N+1}, \quad j = 0, 1, \ldots. \tag{6}$$

***Proof:*** The procedure is quite similar for both numbers, thus we are showing only it for $O_j(N)$. This number can be expressed as

$$\begin{aligned} O_j(N) &= \sum_{r=0}^{j} N^{2r} - \sum_{r=1}^{j} N^{2r-1} = \sum_{r=0}^{j} N^{2r} - N^{-1} \sum_{r=0}^{j} \left(N^{2r-1} - 1\right) \\ &= \frac{1 - N^{2r+2}}{1 - N^2} - N^{-1} \left( \frac{1 - N^{2r+2}}{1 - N^2} - 1 \right) \\ &= \frac{N^{2r+1} + 1}{N+1}. \end{aligned} \tag{7}$$

In order to show that the number is a positive integer it is enough to realize that

$$O_j(N) = N^{-1} + \sum_{r=0}^{j} N^{2r} - N^{-1} \sum_{r=0}^{j} N^{2r-1} > 0, \tag{8}$$

which proves the result. In a similar way the proof for $E_j(N)$ is conducted. □

Now we prove the following result that relates the sequences $S_O(N)$ and $S_E(N)$ with the number of walks in complete graphs. Let us consider a simple graph without multiple links or self-loops $G = (V, E)$ with nodes (vertices) $v_i \in V, i = 1\ldots, n$ and edges $\{v_i, v_j\} \in E$.

A walk of length $k$ is a sequence of (not necessarily distinct) nodes $v_0, v_1, \cdots, v_{k-1}, v_k$ such that for each $i = 1, 2 \cdots, k$ there is an edge from $v_{i-1}$ to $v_i$. If $v_0 = v_k$, the walk is named a *closed* *walk*. A complete graph $K_n$ is the graph having $n$ nodes and every pair of nodes is connected by an edge.

**Theorem 2:** The sequence $S_o(N)$ and $S_E(N)$ give, respectively, the number of walks of odd and even lengths between pairs of nodes in a complete graph $K_{N+1}$.

***Proof:*** Let $M^{2r+1}(p,q)$ be the number of odd walks of length $2r+1$ between the nodes $p$ and $q$ in $K_{N+1}$. Then, it is known that

$$M^{2r+1}(p,q) = \sum_{j=1}^{N+1} \varphi_j(p)\varphi_j(q) \lambda_j^{2r+1}, \qquad (9)$$

where $\varphi_j(p)$ is the $p$th entry of the orthonormalized eigenvector associated with the $\lambda_j$ eigenvalue.

The principal eigenvalue of the adjacency matrix of $K_{N+1}$ is $N$ and its corresponding orthonormalized eigenvector is $\left(\sqrt{N+1}\right)^{-1} \mathbf{1}$. The rest of the eigenvalues are equal to $-1$. Then,

$$\begin{aligned} M^{2r+1}(p,q) &= \frac{1}{\sqrt{N+1}} \frac{1}{\sqrt{N+1}} N^{2r+1} + \sum_{j=2}^{n} \varphi_j(p)\varphi_j(q)(-1)^{2r+1} \\ &= \frac{N^{2r+1}}{N+1} - \sum_{j=2}^{n} \varphi_j(p)\varphi_j(q) \\ &= \frac{N^{2r+1}+1}{N+1}. \end{aligned} \qquad (10)$$

In a similar way the result for the number of walks of even length $M^{2r}(p,q)$ between two nodes in $K_{N+1}$ is proved. □

### 3. Main Result

Let $S_o(N) = O_1(N), O_3(N), \cdots, O_{2r}(N), \cdots$ be the sequence of walks of odd length in a complete graph $K_N$. Thus

$$O_j(N) = \frac{N^{2j+1}+1}{N+1}. \tag{11}$$

Let $G$ be a graph with adjacency matrix $A = (a_{ij})$. We say that $G$ is an *unfolded complete graph* if there exist some $N$, called the *size* of $G$, such that for every edge $p-q$ the virtual power $a_{pq}^{(2j+1)} = O_j(N)$, for every $j \geq 1$. The main result of this section is the following:

**Theorem 1**: Let $G$ be a connected bipartite unfolded complete graph of size $N$. Then $G = K_2 \otimes K_N$.

We divide the proof in several steps.

**Proposition 2**: Let $G$ be a connected, bipartite almost complete graph of size $N$. Then the following hold:

(a) The spectral radius $\rho = N-1$;

(b) $G$ is a regular graph with constant degree $k_p = N-1$ for any node $p$.

**Proposition 3**: Let $G$ be a connected, bipartite graph satisfying:

(c) $G$ is a regular graph with constant degree $k_p = N-1$ for any node $p$;

(d) For every edge $p-q$ we have $a_{pq}^{(2j+1)} = \dfrac{(N-1)^{2j+1}+1}{N}$ for $j = 1, 2$.

Then $G = K_2 \otimes K_N$.

We proceed with the proof of the Propositions.

***Proof of Proposition 2***: Let us write $B = A(K_N) = (b_{rs})$ for the adjacency matrix of $K_N$. Then for any edge $p-q$ in $G$ we have

$$\rho = \lim_{j \to \infty} \sqrt[2j+1]{a_{pq}^{(2j+1)}} = \lim_{j \to \infty} \sqrt[2j+1]{b_{rs}^{(2j+1)}} = N-1. \tag{12}$$

Let $p$ be any vertex of $G$ and choose an edge $p-q$. The matrix $\tilde{A} = \frac{1}{\rho} A$ is double stochastic, its $(p,q)$ entry $\tilde{a}_{pq} = \frac{1}{N-1}$ measures the probability to go from node $p$ to node $q$ in the graph $G$. Therefore $k_p = N-1$. □

We need the following Lemma for the proof of Proposition 2.

**Lemma 4**: Let $G$ be a bipartite graph satisfying (b) and (c) above. Then for any pair of nodes $(p,q)$ in $G$ the following holds:

$$a_{pq}^{(2)} \leq N-2.$$

**Proof**: Assume otherwise that $a_{pq}^{(2)} > N-2$. Observe that $a_{pq}^{(2)} = \sum_s a_{ps} a_{sq} \leq N-1$. Hence $a_{pq}^{(2)} = N-1$ and we get the following full subgraph $H$ of $G$:

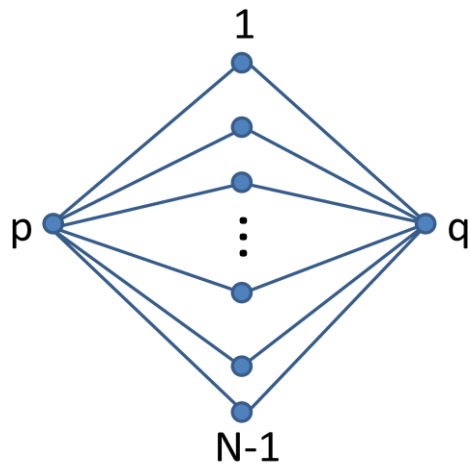

where all neighboring vertices to $p$ and $q$ are $1, \ldots, N-1$. Let us calculate $a_{p1}^{(3)}$. Indeed, using hypothesis (c) we get

$$(N-1)(N-2)+1 = a_{p1}^{(3)} = \sum_{j=2}^{N-1} a_{j1}^{(2)} + N - 1, \tag{13}$$

which yields either one of the following two situations:

1) $a_{j1}^{(2)} = N-1$, for some $1 \neq j$. Without loss of generality, for $j = 2$ we have

$$N-1 = a_{12}^{(2)} = \sum_s a_{1s} a_{s2} \tag{14}$$

Since there are no edges connecting any two of $1,\ldots,N-1$, because $G$ is bipartite, there should exist vertices $2',\ldots,(N-2)'$ not in $H$, maybe not pairwise different, such that $a_{1s} a_{s2} = 1$, for each $s = 2',\ldots,(N-2)'$. We can assume that $a_{1j} a_{j2} = 1$, for each $j = 2,\ldots,N-2$ and there are no paths of length two between $1$ and $N-1$. Repeating the calculation done in (13) we get

$$(N-1)(N-2)+1 = a_{p1}^{(3)} \leq \sum_{j=2}^{N-2} a_{j1}^{(2)} + N - 1 \leq (N-1)(N-3) + N - 1, \tag{15}$$

which is a contradiction.

2) For $i \neq j$ in the set $1,\ldots,N-1$, we have $a_{ij}^{(2)} = N-2$. Then

$$\frac{(N-1)^5 + 1}{N} = a_{p1}^{(5)} \geq \sum_{j=1}^{N-1} a_{pj}^{(3)} a_{j1}^{(2)} \geq \frac{(N-1)^3 + 1}{N} \left[(N-1)(N-2) + N - 1\right], \tag{16}$$

which also yields a contradiction. □

***Proof of Proposition 3***: We shall now describe the structure of the matrix $A^2$. For that purpose we shall calculate all virtual powers $a_{ij}^{(2)}$.

First, for each vertex $p$ we have $a_{pp}^{(2)} = k_p = N-1$. Let $p,q$ be different vertices in $G$. By the Lemma above, $a_{pq}^{(2)} \leq N-2$. Assume that $a_{pq}^{(2)} \neq 0$ and take any vertex $j$ with edges $p-j-q$. Then for the neighbours $1,\ldots,N-1$ of $j$ we get, using Proposition 2,

$$(N-1)(N-2)+1 = a_{pj}^{(3)} \leq \sum_{p \neq j=1}^{N-1} a_{pi}^{(2)} + N - 1 \leq (N-2)^2 + N - 1, \tag{17}$$

which implies equality $a_{pq}^{(2)} = N-2$.

We observe that the number of vertices of the graph $G$ is $n = 2(N-1)$. Indeed, since $G$ is bipartite, we get two non-empty classes $V_1$, $V_2$ of totally disconnected vertices in the set $1,\ldots,n$. Let $p, q, s$ be three different vertices in $V_1$ and $p-j-q-i$ edges. Then

$$a^{(2)}_{pq} = N - 2 = a^{(2)}_{qs}. \tag{18}$$

Since this number is not zero, there is a vertex $k$ in $V_2$ which is connected to $p$ and $s$, and $a^{(2)}_{ps} \neq 0$. Therefore $A^2$ has $a^{(2)}_{pq} = N - 2$ or 0 depending if $p$, $q$ belong or not to the same class $V_i$. Hence each of $V_1$ and $V_2$ have the same cardinality $N-1$.

Moreover it is clear that the adjacency matrix $A$ has entries $a_{pq} = 0$ or 1 depending if $p$, $q$ are different and in distinct classes or belong to the same class $V_i$. As desired, it follows that $A = A(K_2) \otimes A(K_n)$. □

## 4. Graph Unfolding Construction

Despite it is straightforward to obtain the complete unfolded graphs from the Kronecker product indicated in the previous section we design a method here that allows to 'visualizing' and give more intuition about the transformation described previously. The method developed here has some resemblances with the edge subdivision of a graph. The subdivision graph $S(G)$ of a graph $G$ is the graph obtained subdividing every edge of $G$. The subdivision of the edge $\{p, q\}$ produces a graph containing one new node $r$, and the edge $\{p, q\}$ is replaced by $\{p, r\}$ and $\{r, q\}$. Here we consider a complete graph $K_N$ and select a complete subgraph $S_{N-1}$ having $N-1$ nodes. We then add a new node $r$ and replace all the existing edges in $S_{N-1}$ by new edges linking the node $r$ with every node of $S_{N-1}$. This procedure is repeated for every complete subgraph $S_{N-1}$ having $N-1$ nodes in $K_N$ (see

Figure 1). The subdivison is the particular case for $K_3$ in which every subgraph $S_2$ is an edge which is subdivided by new nodes. The general algorithm is described below:

*Complete Graph Unfolding Algorithm*

---

1) Consider a complete graph $K_N$,

2) Select a subgraph $S$ having $N-1$ nodes in $K_N$,

3) Insert a new node which is attached to every node of the subgraph $S$,

4) Repeat 2) and 3) for all the distinct subgraphs $S$ having $N-1$ nodes in $K_N$,

5) Remove all the links of $K_N$. The resulting graph is $G_{2N}$.

---

**Insert Figure 1 about here**

By using the complete graph unfolding algorithm we can 'visualize' the graph operations involved. For instance, for $N=3$ two triangles $K_3$ give rise to a hexagon and for $N=4$ two tetrahedra give rise to a cube (see Figure 2). For $N \geq 5$ the complete graphs can only embedded into high-dimensional spaces, which make impossible a simple visualization as before.

**Insert Figure 2 about here.**

**5. Characteristics of unfolded complete graphs**

The graph $G_{2N}$ obtained (see some examples in Figure 1) from this complete graph unfolding transformation has the following general characteristics:

1) $G_{2N}$ has $2N$ nodes, where $N$ is the number of nodes in the original $K_N$;

2) $G_{2N}$ is regular with degree $k = N-1$;

3) $G_{2N}$ is bipartite;

4) the adjacency matrix of $G_{2N}$ is

$$A(G_{2N}) = \begin{pmatrix} 0 & A(K_N) \\ A(K_N) & 0 \end{pmatrix}. \tag{19}$$

**Lemma 5:** The distance matrix of the graph $G_{2N}$ has the following block form

$$D(G_{2N}) = \begin{pmatrix} 2A(K_N) & 3I + A(K_N) \\ 3I + A(K_N) & 2A(K_N) \end{pmatrix}. \tag{20}$$

***Proof:*** The graph $G_{2N}$ is formed by two disjoint sets $V_1$ and $V_2$ of $N$ nodes each. Every node in $V_1$ is connected to $N-1$ nodes in $V_2$. Then, any pair of nodes in $V_1$, respectively in $V_2$, are separated at distance two. This gives rise to the terms $2A(K_N)$ in the main diagonal of $D(G_{2N})$. Let us arrange the labels of the nodes in $G_{2N}$ as follows. We label the nodes in $V_1$ by $1, 2, \cdots, N-1$ and those in $V_2$ by $1', 2', \cdots, N'-1$. Then, we connected every node in $V_1$ to every node in $V_2$ except for the pairs $v, v'$. Consequently, every pair $v, u'$ is separated at distance one (they are connected), while the pairs $v, v'$ are at distance three. Then, the term $3I + A(K_N)$ appears as the nondiagonal entries of $D(G_{2N})$, which proves the result. □

**Corollary 6:** The diameter of the graphs $G_{2N}$ is 3.

**Lemma 7:** The number of walks in the graph $G_{2N}$ is given by:

$$M_{ij}^r = \begin{cases} \dfrac{k^{2r+1}+1}{k+1}, & \text{if } i \neq j,\ d_{ij} = 1 \text{ and } r \text{ is odd} \\[2mm] \dfrac{k^{2r}+1}{k+1}, & \text{if } i \neq j,\ d_{ij} = 2 \text{ and } r \text{ is even} \\[2mm] 0, & \text{if } i = j,\ \text{and } r \text{ is odd} \\[2mm] \dfrac{k^{2r}+1}{k+1} + 1, & \text{if } i = j,\ \text{and } r \text{ is even} \end{cases}$$

where $k = \dfrac{N}{2} - 1$ is the degree of a node in $G_{2N}$.

***Proof***: The number of walks of odd length between pairs of adjacent nodes in $G_{2N}$ was already proved in Theorem 1. In order to find the number of walks of even length we consider the block form of $A(G_{2N})$, such that

$$A(G_{2N})^{2r+2} = \begin{pmatrix} A(K_N)^{2r+2} & 0 \\ 0 & A(K_N)^{2r+2} \end{pmatrix}. \tag{21}$$

Using Lemma 4 we can see that the nodes for which $A(K_N)^{2r+2} \neq 0$ are those separated at distance two. Consequently, the number of walks of even lengths between pairs of nodes at distance two in $G_{2N}$ is the same as the number of even walks between two nodes in $K_N$. Using Theorem 2 the result is completed.

The final statement is proved as follows. First, we can see that the number of closed walks of even length is zero as resulting from the main diagonal of

$$A(G_{2N})^{2r+1} = \begin{pmatrix} 0 & A(K_N)^{2r+1} \\ A(K_N)^{2r+1} & 0 \end{pmatrix}. \tag{22}$$

From (21) it can be seen that $\left[A(G_{2N})^{2r+2}\right]_{ii} = \left[A(K_N)^{2r+2}\right]_{ii}$. Then,

$$\begin{aligned} \left[A(K_N)^{2r+2}\right]_{ii} &= \frac{1}{N}(N-1)^{2r+2} + \sum_{j=2}^{N}\left[\varphi_j(i)\right]^2 \\ &= \frac{1}{N}(N-1)^{2r+2} + 1 - \frac{1}{N} \\ &= \frac{(N-1)^{2r+2} - 1}{N} + 1. \quad \square \end{aligned} \tag{23}$$

Now we will prove that the graphs $G_{2N}$ are the incidence matrix of certain symmetric block designs. A block design [4-6] with parameters $(v, k, \lambda)$, also called a $2-(v, k, \lambda)$ design, is a finite point set of cardinality $v$, a collection of blocks, which are subsets of the point set and have the following properties: i) each block has cardinality $k$, ii) each (unordered) pair of points appears in exactly $\lambda$ blocks. Block designs have important applications in different areas such as experimental design, software testing, cryptograph, algebraic geometry, among others [5].

**Theorem 8:** The graphs $G_{2N}$ are the incidence graphs of symmetric 2-$(k+1, k, k-1)$ designs, where $k = \dfrac{N}{2} - 1$ is the degree of a node in $G_{2N}$ and $N$ is the number of nodes in $k = \dfrac{N}{2} - 1$ is the degree of a node in $K_N$.

***Proof:*** First we find the spectra of $G_{2N}$. If the eigenvalues of $P$ and $Q$ are $\lambda_1, \cdots, \lambda_n$ and $\mu_1, \cdots, \mu_m$, respectively, then the eigenvalues of $P \otimes Q$ are $\lambda_i \mu_j$, $i = 1, \cdots, n, j = 1, \cdots, m$. Consequently, the eigenvalues of $G_{2N}$ are $sp(G_{2N}) = \{1 \cdot [k]^1, 1 \cdot [-1]^{N-1}, -1 \cdot [-1]^{N-1}, -1 \cdot [k]^1\}$, which gives the result. Consequently, the spectra of $G_{2N}$ is

$$sp(G_{2N}) = \{[k]^1, [1]^{N-1}, [-1]^{N-1}, [-k]^1\}. \tag{24}$$

Now, according to a result due to Cvetković, Doob, Sach [9] a connected bipartite graph with four distinct eigenvalues must be the incidence graph of symmetric 2-$(v, r, s)$ designs. As a consequence, the graphs $G_{2N}$ are the incidence graphs of symmetric 2-$(k+1, k, k-1)$ designs, where $k = \dfrac{N}{2} - 1$ is the degree of a node in $G_{2N}$ and $N$ is the number of nodes in $k = \dfrac{N}{2} - 1$ is the degree of a node in $K_N$. □

A graph is distance-regular when, for any two nodes $p$ and $q$ at distance $d_{pq} = l$ the number $c_l$, $a_l$ and $b_l$ of nodes which are adjacent to $p$, and at distance $l-1$, $l$, and $l+1$, respectively, from $q$ only depend on $l$. Fiol [10] has proved that a graph is distance-regular if and only if, for each nonnegative integer $s$, the number $a_{pq}^{(s)}$ of walks of length $s$ between nodes $p$ and $q$ only depends on $d_{pq} = l$.

**Theorem 9:** The graphs $G_{2N}$ are *distance-regular graphs*.

*Proof*: We use the previously mentioned Theorem by Fiol [10] in order to prove this result. From Lemma 7 we know that the number of closed walks is the same for any node. Also from Lemma 7 we know that the number of walks of even length is equal for any pair of nodes at distance two. From Lemma 5 we have that the distance matrix of $G_{2N}$ is:

$$D(G_{2N}) = \begin{pmatrix} 2A(K_N) & 3I + A(K_N) \\ 3I + A(K_N) & 2A(K_N) \end{pmatrix}.$$

Then, because

$$A(G_{2N})^{2r+1} = \begin{pmatrix} 0 & A(K_N)^{2r+1} \\ A(K_N)^{2r+1} & 0 \end{pmatrix},$$

the number of walks of odd length is equal to $\frac{k^{2r+1}+1}{k+1}$ if $d_{pq} = 1$ or to $M_{pp}^{2r+1}(K_N)$ if $d_{pq} = 3$.

Thus, the number $a_{pq}^{(s)}$ of walks of length $s$ between nodes $p$ and $q$ only depends on $d_{pq} = l$.

□

Finally, we show that the graphs $G_{2N}$ have another important property, namely they are Ramanujan graphs. Let $\lambda_j \neq \pm k$ be referred as the nontrivial eigenvalues of a regular bipartite graph $G$. Let us denote by $\lambda(G)$ the maximum of the absolute value of all the nontrivial eigenvalues of $G$. Then, a Ramanujan graph is a $k$-regular graph satisfying [8]: $\lambda(G) \leq 2\sqrt{k-1}$. Then, we have the following result. From Theorem 8 we can see that $\lambda(G_{2N}) = 1$ and $1 \leq 2\sqrt{k-1}$ for any $k \geq 2$. Thus, we can write this result as follow.

**Lemma 10**: The graphs $G_{2N}$ are Ramanujan graphs for $k \geq 2$.

More importantly, the graphs $G_{2N}$ have very large edge expansion constants [11, 12]. For a given graph $G = (V, E)$ and for two disjoint subsets $P$ and $Q$, let us denote by

$e(P,Q)$ the number of edges between $P$ and $Q$. The edge expansion constant $i(G)$ is defined by

$$i(G) = \min \frac{e(P, V-P)}{|P|},$$

where the minimum is taken over all non-empty subsets $P$ and $Q$ of cardinality at most $|V|/2$. The celebrated Alon-Milman theorem [13] proves that

$$i(G) \geq \frac{k - \lambda(G)}{2}.$$

Thus, for $G_{2N}$ the expansion constant is bounded as

$$i(G_{2N}) \geq \frac{N-4}{4}.$$

## 6. Closing remark

We have introduced two new integer sequences which we prove correspond to the sequences of even and odd walks in complete graphs $K_n$. The fundamental question posted here is whether there are other graphs with exactly the same sequence of walks as those of the complete graphs. We have proved that the answer is positive and the graphs having such sequences are unique. They can be constructed as the Kronecker product $G_{2n} = K_2 \otimes K_n$. Such graphs are the incidence graphs of block designs, are distance-regular and Ramanujan graphs, which open a wide horizon of applications in many areas of research.

**Acknowledgement**

EE thanks CIMAT for warm hospitality during January-February 2013 and for a visiting professorship position during this period.

**Figure captions**

Figure 1. Illustration of the graph unfolding process for $N=4$ (top) and $N=5$ (bottom), see text for description.

Figure 2. Geometric interpretation of the result of the graph unfolding process for $N=4$, where two rotated tetrahedra gives rise to a cube.

**Figure 1**

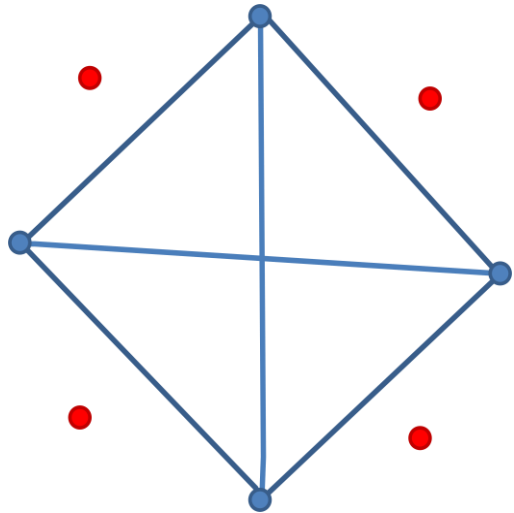 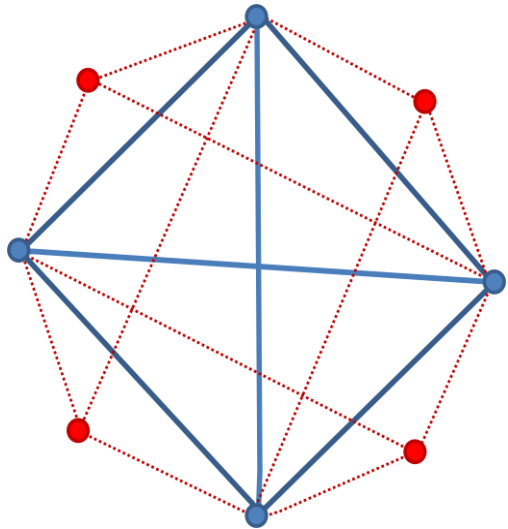 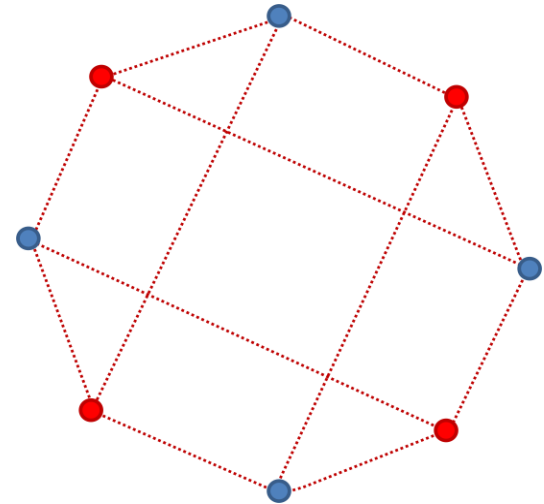

**Figure 2**

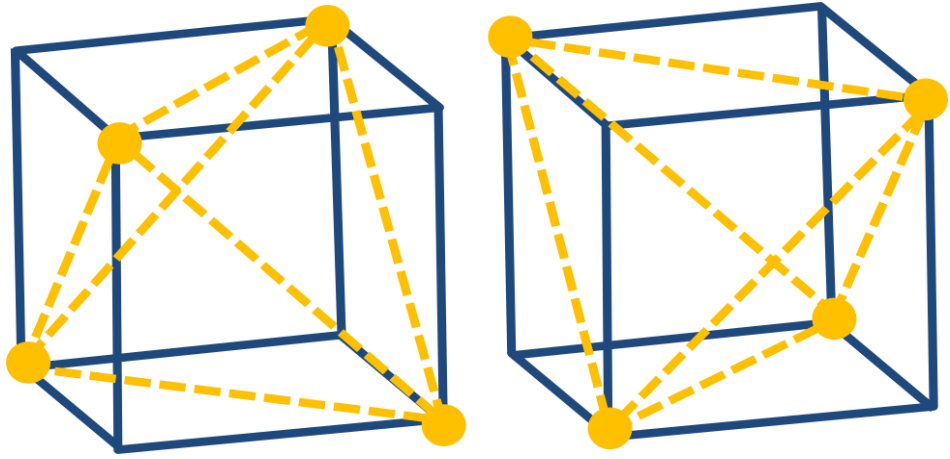